\documentclass{article}
\usepackage[a4paper, total={6in, 8in}]{geometry}
\usepackage{graphicx,subfigure,microtype,wrapfig,natbib,booktabs,hyperref,xcolor}  
\usepackage{amsmath, amsthm, amssymb, mathtools, bbm, bm, upgreek}
\usepackage{algorithm,algorithmic}
\usepackage{caption,subcaption, enumitem}
\usepackage[capitalize,noabbrev]{cleveref}
\usepackage[textsize=tiny]{todonotes}
\theoremstyle{remboldstyle}
\newtheorem*{definition*}{Defintion}
\newtheorem*{thm*}{Theorem}
\newtheorem*{theorem*}{Theorem}
\newtheorem*{assu*}{Assumption}
\newtheorem*{assumption*}{Assumption}
\newtheorem*{lemma*}{Lemma}
\newtheorem*{example*}{Example}
\newtheorem*{ex*}{Example}
\newtheorem*{prop*}{Proposition}
\newtheorem*{proposition*}{Proposition}
\newtheorem*{coro*}{Corollary}
\newtheorem{observation}{Observation}

\theoremstyle{plain}
\newtheorem{theorem}{Theorem}[section]
\newtheorem{proposition}[theorem]{Proposition}

\newtheorem{lemma}[theorem]{Lemma}
\newtheorem{corollary}[theorem]{Corollary}
\theoremstyle{definition}
\newtheorem{definition}[theorem]{Definition}
\newtheorem{assumption}[theorem]{Assumption}
\newtheorem{remark}[theorem]{Remark}

\newcommand{\su}[1]{}
\newcommand{\new}{\color{black}}

\newcommand{\N}{{\cal N}}
\newcommand{\blue}{\color{blue}}

\newcommand{\bitem}{\begin{itemize}}
\newcommand{\eitem}{\end{itemize}}
\newcommand{\benum}{\begin{enumerate}}
\newcommand{\eenum}{\end{enumerate}}
\newcommand{\blem}{\begin{lemma}}
\newcommand{\elem}{\begin{lemma}}
\newcommand{\bdefn}{\begin{definition}}
\newcommand{\edefn}{\end{definition}}
\newcommand{\bprop}{\begin{proposition}}
\newcommand{\eprop}{\end{proposition}}
\newcommand{\bobsv}{\begin{observation}}
\newcommand{\eobsv}{\end{observation}}

\newcommand{\bcn}{\begin{coro*}}
\newcommand{\ecn}{\end{coro*}}
\newcommand{\bdn}{\begin{definition*}}
\newcommand{\edn}{\end{definition*}}
\newcommand{\bln}{\begin{lemma*}}
\newcommand{\eln}{\end{lemma*}}
\newcommand{\btn}{\begin{theorem*}}
\newcommand{\etn}{\end{theorem*}}
\newcommand{\bpn}{\begin{proposition*}}
\newcommand{\epn}{\end{proposition*}}

\newcommand{\ps}{\begin{proof}[Sketch]}
\newcommand{\brmk}{\begin{remark}}
\newcommand{\ermk}{\end{remark}}
\newcommand{\bcoro}{\begin{corollary}}
\newcommand{\ecoro}{\end{corollary}}
\newcommand{\bcom}{}

\newcommand{\avg}{average}

\newcommand{\assu}{assumption}

\newcommand{\bs}{\backslash}

\newcommand{\cvx}{convex}

\newcommand{\cond}{condition}

\newcommand{\corres}{corresponding}

\newcommand{\distr}{distribution}

\newcommand{\func}{function}

\newcommand{\gs}{\gtrsim}

\newcommand{\ho}{\mathbb}

\newcommand{\indep}{independent}

 \newcommand{\ind}{\perp\!\!\!\!\perp}

\newcommand{\IOW}{In other words}

\newcommand{\ineq}{inequality}

\newcommand{\ls}{\lesssim}

\newcommand{\msb}{measurable}

\newcommand{\mtx}{matrix}

\newcommand{\nbhd}{neighborhood}

\newcommand{\Omg}{\Omega}

\newcommand{\OTOH}{On the other hand}

\newcommand{\rb}{\right}

\newcommand{\lb}{\left}

\newcommand{\Prb}{Probability}
\newcommand{\prb}{probability}

\newcommand{\parti}{particular}

\newcommand{\pg}{\paragraph}

\newcommand{\rv}{random variable}

\newcommand{\rar}{\rightarrow}

\newcommand{\real}{\mathbb{R}}
\newcommand{\resp}{respectively}

\newcommand{\sat}{satisfy}

\newcommand{\satd}{satisfied}

\newcommand{\sse}{\subseteq}
\newcommand{\sps}{suppose}
\newcommand{\Sps}{Suppose}

\newcommand{\strfwd}{straightforward}

\newcommand{\unif}{uniform}
\newcommand{\unk}{unknown}

\newcommand{\vtxs}{vertices}

\newcommand{\wh}{\widehat}

\let\eps\varepsilon

\title{Experimentation Under Non-stationary Interference}
\author{
  Su Jia\thanks{Supply Chain Optimization Technologies (SCOT), Amazon; \texttt{sujiaxxx@amazon.com}} \and
  Peter Frazier\thanks{School of Operations Research and Information Engineering (ORIE), Cornell University} \and
  Nathan Kallus\footnotemark[2] \and
  Christina Lee Yu\footnotemark[2]
}

\date{}

\begin{document}
\maketitle

\begin{abstract} 
We study the estimation of the {\em average treatment effect} (ATE) in randomized controlled trials under a dynamically evolving interference structure.
This setting arises in applications such as ride-sharing, where drivers move over time, and social networks, where connections continuously form and dissolve.
In particular, we focus on scenarios where outcomes exhibit spatio-temporal interference driven by a sequence of random interference graphs that evolve independently of the treatment assignment.
Loosely, our main result states that a truncated Horvitz-Thompson estimator achieves an MSE {\new that vanishes linearly in the number of spatial and time blocks, times a factor that measures the average complexity of the interference graphs.
As a key technical contribution (that contrasts the static setting \citealt{jia2023clustered}), we present a fine-grained covariance bound for each pair of space-time points that decays exponentially with the time elapsed since their last ``interaction''. 
Our results can be applied to many concrete settings and lead to simplified bounds, including where the interference graphs (i) are induced by moving points in a metric space, or (ii) follow a dynamic Erdos-Renyi model, where each edge is created or removed \indep ly in each time period.
\su{... And the ``complexity" terms shouldn't depend on N,T. 
Also, state the results in a general way so that we can say something like ``for specific interference graph we have xx MSE''. 
E.g. examples where we should choose non-uniform spatial partitioning.  
Need to make more interpretable.
Bottom-up: use interpretable special case to attract readers.}
}
\end{abstract}

\section{Introduction}
Experimentation under interference — that is, when the {\em Stable Unit Treatment Value Assumption} (SUTVA) \citep{rubin1978bayesian} is violated — poses a fundamental challenge in modern experimentation:
In many large-scale online environments such as marketplaces, the treatment assigned to one unit may affect the outcomes of others through shared resources, social connections, or platform algorithms. 
Ignoring such spillover effects can lead to biased estimates, invalid inference, and misguided product decisions. 
As a result, modeling and mitigating interference has become a central challenge in the design of trustworthy experimentation pipelines.

Most existing work in this area focuses on either (1) a single-period setting where non-stationarity is irrelevant, or (2) a multi-period model with a {\bf static} interference pattern, usually characterized by a graph.
However, such assumptions can be restrictive for modern applications, where user behavior and network structure evolve over time. 
In real-world systems, interference can evolve as users enter or exit the platform, content recommendations adapt, or interactions become more or less frequent. 
These dynamic patterns make it essential to develop experimental methods that account for both temporal evolution and the fluid nature of the interference graph.

We initiate a study of experimentation under a {\em non-stationary} interference. 
Non-stationarity is {\bf two-fold}: We allow not only the outcome distributions to vary (across  both space and time), but also the interference graph to evolve over time, {\bf independently} of the treatments.
More precisely, conditional on the interference graphs, the outcomes of  each fixed individual follow a Markov reward process, whose state transitions may depend on the treatment of other individuals in the same time period.

\subsection{Our Contributions}
\su{expand}
We progressively show the following results on the Horvitz-Thompson (HT) estimator.
\benum 
\item {\bf Covariance Bound for the HT Estimators.}  We establish a general covariance bound in terms of a novel quantity called the \emph{last interaction time} (LIT). Specifically, we show that the covariance between the HT estimators of individual treatment effects decays exponentially with the time elapsed since their LIT (\cref{prop:cov_LIT}). 
This is a key technical component for our main result and is of independent interest, as it holds for {\bf any} experimental design (i.e. distribution over the treatment assignment matrix).
\item {\bf Bounding MSE Using the Average Cluster Degree.} We present a general MSE bound for a broad class of designs, which we call {\em vertical designs}: Partition the time horizon uniformly into time blocks, then within each time block, partition the population into spatial blocks and assign treatment/control at the block level.
We show that the MSE of the HT estimator is  $1/NT$ times the {\em average cluster degree}, where the average is taken over all individuals and time blocks (\cref{thm:mse}).
Here, two individuals are called {\em cluster-neighbors} in a time block if their neighborhoods intersect the same spatial block, and the cluster degree of an individual is the number of cluster-neighbors.
\eenum 

\subsection{Related Work}
\su{This lit rev is copied from our earlier paper. Will make minor edits.}

\noindent{\bf Violation of SUTVA.}
Experimentation is broadly deployed in e-commerce that is simple to execute \citep{kohavi2017surprising,thomke2020building,larsen2023statistical}.
As a key challenge, the violation of the {\em Stable Unit Treatment Value Assumption} (SUTVA) has been viewed as problematic in online platforms \citep{blake2014marketplace}. 
Many existing works that tackle this problem assume that interference is summarized by a low-dimensional exposure mapping and that individuals are individually randomized into treatment or control by Bernoulli randomization \citep{manski2013identification,toulis2013estimation,aronow2017estimating,basse2019randomization,forastiere2021identification}.
Some work departed from unit-level randomization and introduced cluster dependence in unit-level assignments in order to improve estimator precision, including 
\citealt{ugander2013graph,jagadeesan2020designs,leung2022rate,leung2023network}, just to name a few.

There is another line of work that considers the temporal interference 
(or {\em carryover} effect).
Some works consider a fixed bound on the persistence of temporal interference (e.g., \citet{bojinov2023design}), while other works considered temporal interference arising from the Markovian evolution of states
\citep{glynn2020adaptive,farias2022markovian,johari2022experimental,shi2023dynamic,hu2022switchback,li2022network,li2023experimenting}.
Apart from being limited to the single-individual setting, many of these works differ from ours either by (i) focusing on alternative objectives, such as stationary outcome \citet{glynn2020adaptive}, or (ii) imposing additional assumptions, like observability of the states \citet{farias2022markovian}. 

\

\noindent{\bf Spatio-temporal Interference.}
Although extensively studied separately and recognized for its practical significance, experimentation under spatio-temporal interference has received relatively limited attention  previously.
Recently, \citet{ni2023design} attempted to address this problem, but their carryover effect is confined to one period.
 Another closely related work is \citealt{li2022network}. 
Similar to our work, they specified the spatial interference using an interference graph and modeled temporal interference by assigning an MDP to each individual. 
In our model, the transition probability depends on the states of all neighbors (in the interference graph).  
In contrast, their evolution  depends on the {\bf sum} of the outcome of direct neighbors. 
Moreover, our work focuses on ATE estimation for a fixed, unknown environment, whereas they focus on the large sample asymptotics and mean-field properties.
\citet{wang2021causal} studied direct treatment effects for panel data under spatiotemporal interference, but focused on {\em asymptotic} properties instead of finite-sample bounds.

\

\noindent{\bf Off-Policy Evaluation (OPE)}
Since we model temporal interference using an MDP, our work is naturally related to reinforcement learning (RL).
In fact, our result on ATE estimation can be rephrased as OPE 
\citep{jiang2016doubly,thomas2016data} in a multi-agent MDP: Given a {\em behavioral} policy from which the data is generated, we aim to evaluate the mean reward of a {\em target} policy. 
The ATE in our work is essentially the difference in the mean reward between two target policies  
(all-1 and all-0 policies), and the behavioral policy is given by clustered  randomization.
However, these works usually require certain states to be observable, which is not needed in our work.
Moreover, these works usually impose certain assumptions on the non-stationarity, which we allow to be completely arbitrary.
Finally, we focus on rather general data-generating policies (beyond fixed-treatment policies) and estimands (beyond ATE), compromising the strengths of the results.

\su{Add these: Decaying interference: 
\cite{faridani2024} }




\section{Formulation}
\label{subsec:basic_setting}
We consider a multi-period model that captures both spatial interference and temporal interference. 
Our framework naturally extends that of \cite{hu2022switchback,jia2023clustered}.
\begin{itemize}[noitemsep, topsep=0pt]
\item {\bf Basics.}  Consider $N$ individuals, e.g., drivers on a ride-sharing platform, or users in an online platform, and $T$ rounds (e.g., days). 
Denote the two {\em treatment arms} (or, simply, {\em arms}) by $0,1$.\footnote{We consider two arms for simplicity. 
It is \strfwd\ to extend to multiple arms.}
A {\em treatment assignment} is a \mtx\ $W\in \{0,1\}^{N\times T}$.
\item {\bf Spatial Interference via Interference Graphs.}
Interference between individuals is governed by a sequence of {\em interference graphs} $G_t = ([N],\ E_t)$.\footnote{Every node always has a self-loop (i.e., an edge connecting the node to itself).}
For now, we assume that each $G_t$ is {\bf deterministic}, and explain later why it is \strfwd\ to extend our results to the case where $G_t$'s evolve randomly (\indep ly of $W$).
Denote $\N_t(i) = \{i\}\cup \{j\in [N]: (i,j)\in E_t\}$ the neighbors of $i$ at time $t$.
\item {\bf Temporal Interference via MDPs.}
Each individual has a (possibly unobservable) {\em state} $S_{it}\in {\cal S}$ where $\cal S$ is a (possibly \unk) state space.
The state evolves in a Markovian fashion. 
More precisely, $S_{i,t+1}$ depends only on (i) $S_{it}$ and (ii) the treatments of $i$'s neighbors at time $t$, i.e., $W_{\N_t(i)}$ but not on the state or treatment prior to $t-1$.
Formally, conditional on the {\em treatment assignment} $w \in \{0,1\}^N$ for round $t$, the state transition of individual $i$ is governed by an \unk\ {\em transition kernel} 
\[
P^w_{it}(s,\cdot) := P_{it}\lb[S_{i,t+1}\in \cdot \mid S_{it}=s,\ W_{{\cal N}_t(i), t} = w_{{\cal N}_t(i), t} \rb].\footnote{
We use lower and upper case for    deterministic and random variables \resp.}
\]
\item {\bf Outcomes.} At the end of a time period,  we observe a random {\em outcome}
\begin{align}\label{eqn:010724}
Y_{it} = \mu_{it}\lb(S_{it}, w_{{\cal N}_t(i),t}\rb) + \eps_t 
\end{align}
where $\mu_{it}: {\cal S} \times \{0,1\}^{ {\cal N}_t(i)} \rar [0,1]$ is an \unk\ {\em outcome \func} and $\eps_t$ is an ``external'' error \indep\ of any other variables, with mean $0$ and variance at most $\sigma^2$.
\end{itemize}
Denote by $\Delta(\cdot)$ the \prb\ simplex and $d_{\rm TV}$ the TV-distance.
As in \citealt{hu2022switchback},   assume:

\begin{assumption}[Rapid-mixing] 
\label{assu:rapid_mixing}
There is a constant $t_{\rm mix}>0$ s.t. for any $i,t$, $f,f'\in \Delta({\cal S})$, and $w\in \{0,1\}^N$, we have
\[
d_{\rm TV}(f P^w_{it}, f' P^w_{it})= \frac 12 \|f P^w_{it} - f' P^w_{it}\|_1 
\le e^{-1/t_{\rm mix}} \cdot \|f - f'\|_1.
\]
\end{assumption}

We aim to estimate the difference in the \avg\ outcomes between all-treatment vs. all-control.

\bdefn[Average Treatment Effect]
For any $f \in \Delta({\cal S})$, we define 
\[{\rm ATE} (f):= \ho{E}\lb[\frac 1{NT} \sum_{i,t} Y_{it}\ \Big| \  W = {\bf 1}_{N\times T},\ S_0 \sim f\rb]- \ho{E}\lb[\frac 1{NT} \sum_{i,t} Y_{it}\ \Big| \  W= {\bf 0}_{N\times T},\ S_0 \sim f\rb].\]
\edefn 

Thanks to the rapid-mixing \assu, 
it is not hard to show that the initial state \distr\ $f$ only matters up to a lower-order term in the MSE.
We will therefore suppress the initial \distr.

\bprop[Initial State Doesn't Matter Much] \Sps\ \cref{assu:rapid_mixing} holds. Then, for any $f,f' \in \Delta(\cal S)$, we have 
\[|{\rm ATE}(f) - {\rm ATE}(f')| \le \frac 1{NT} \sum_{t=1}^T e^{-t/t_{\rm mix}} \ls \frac {t_{\rm mix}}{NT}.\]
\eprop


\section{Preliminaries}
To state our MSE bound, we need to introduce several key components, including standard ones such as randomized design and Horvitz-Thompson estimator, and a novel concept called the maximum \indep\ time. 

\subsection{Design and Estimator} 
A {\em randomized design} (or  simply {\em design}) is a \prb\ distribution over $\{0,1\}^{N \times T}$.  
Consider a natural class of designs: First partition time into time block. 
At the {\bf start} of each time block, partition the individuals into  spatial blocks.\footnote{
The interference graphs can be completely arbitrary and need not be embedded in any underlying metric space. 
Accordingly, ``spatial blocks'' are not spatial in any geometric sense — they are simply subsets of individuals.}
Then, randomly assign arms at the spatial block. Each individual will then carry this assigned arm {\bf throughout} the time block.

\bdefn[Vertical Design] 
For each $k=1,\dots, T/\ell$, let $\Pi_k$ be a partition of $[N]$ into {\em spatial blocks}. 
For each time block $k$ and spatial block $B\in \Pi_k$, we call $B\times [(k-1)\ell+1,\ k\ell]$ a {\em spatio-temporal cluster} (or {\em cluster}).
We \indep ly assign each cluster $C$ an arm $Z_C\sim {\rm Ber}(\frac 12)$.\footnote{We chose $1/2$ for simplicity. 
It is \strfwd\ to extend to other Bernoulli means.}
Then, set $W_{it} = Z_{C[i,t]}$ for each $i,t$, where $C[i,t]\sse [N]$ is the (unique) cluster containing   $i$ at time $t$.
\edefn



\bdefn[Spatio-temporal Neighborhood] 
For any $i,t$ and $r>0$, we define the {\em spatio-temporal neighborhood} as 
\[
{\N}^r(it) := \bigcup_{\tau = t-r}^t \lb(\N_\tau(i)\times \{\tau\}\rb).
\]
\edefn

We will show that the following estimator achieves an MSE that decays  linearly in $N,T$ times a factor that depends on the \avg\  complexity of the interference graphs.
This estimator is based on an exposure mapping that is {\bf mis-specified}, in that we truncate the history prior to round $t-r$.

\bdefn[Exposure Mapping]
For any $i,t,a$ and $r>0$, we define the {\em exposure mapping} as\footnote{Be careful: Note that $p_{it}^r$ is {\bf deterministic}, so we can factor them out in the covariance.
If the interference graph changes {\bf randomly}, then $p_{it}^r$ is a \rv\ as it depends on $W$. 
In this case, ${\rm Cov}(\wh \Delta_{it},\wh \Delta_{i't'})$ may not decay exponentially in $|t-t'|$, since $p_{it}^r$ and $p_{i't'}^r$ may have a high covariance even if $|t-t'|$ is huge.
This can be addressed, though, by choosing a larger but deterministic $p_{it}^r$.}
\[
X_{ita}^r(w) := \mathbbm{1} (w_{\N^r(it)} \equiv a\cdot {\bf 1})
\] 
for each $w\in \{0,1\}^{N\times T}$. 
For a design $\cal D$, we define the {\em exposure \prb} as 
\[
p_{ita}^r := \ho{P}_{W\sim \cal D}\lb[X_{ita}^r(W) =1\rb].
\]
\edefn

Note that since we assumed the \prb\ both arms are $0.5$, we have $p_{it1}^r = p_{it0}^r =: p_{it}^r$.

\bdefn[Horvitz-Thompson Estimator]
For each $i,t,a$, we define 
\[
\wh Y_{ita}^r = \frac{X_{ita}^r}{p_{it}^r} Y_{it}.
\]
The {\em Horvitz-Thompson} (HT) estimator is
\[\wh \Delta^r =\frac 1{NT}\sum_{(i,t)\in [N]\times [T]} \wh \Delta^r_{it} \quad {\rm where} \quad \wh \Delta_{it}^r = \wh Y_{it1}^r - \wh Y_{it0}^r.\]
\edefn
We will refer to these $\wh \Delta_{it}^r$ as the {\em HT terms}.
We will show that regardless  of the design, the HT estimator achieves low bias for an appropriate choice of $r$.

\begin{proposition}[Bias bound] \label{prop:bias} 
For any $r\ge 0$ and design $\cal D$, we have 
\[\ho{E}\lb[\lb|\wh \Delta^r - \Delta\rb|\rb] \le e^{-r/t_{\rm mix}}.\]
\end{proposition}
The  proof is \strfwd; we defer it to \cref{sec:bias}.

\subsection{Clustering-Induced Graphs and Cluster Degree}

\begin{figure}
\centering    \includegraphics[width=0.7\textwidth]
{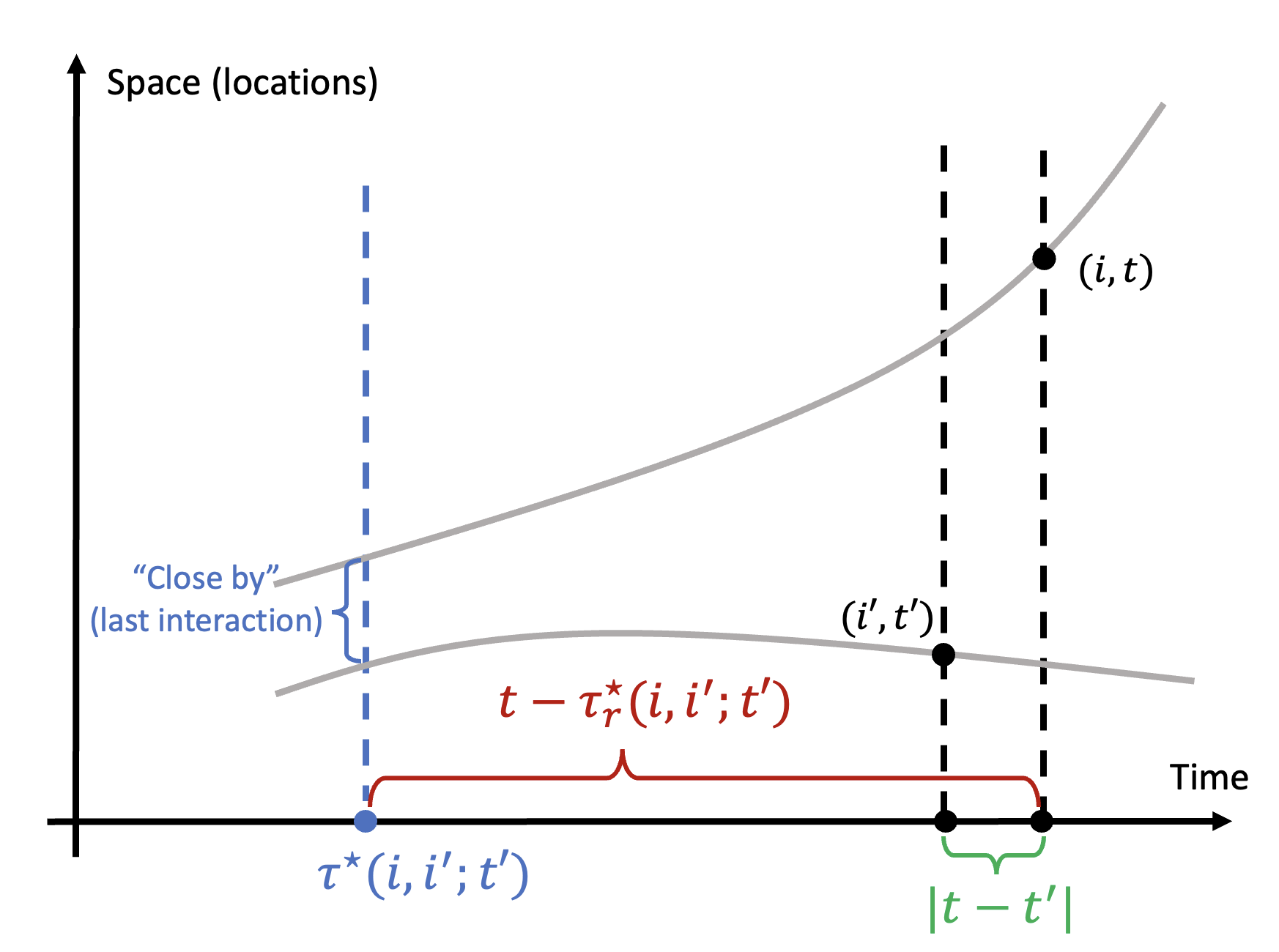}
\caption{{\bf Naive Covariance Bound Is Loose.} 
Wlog assume that $t'\le t$. The naive covariance bound only gives $e^{-|t-t'|/t_{\rm mix}}$ (see \cref{eqn:111324}), which can be loose if $t\approx t'$. 
To tighten this, we will define the last interaction time $\uptau^\star(i,i';t')$ and show that the covariance decays exponentially in $t -\uptau^\star(i,i';t')$. }
\label{fig:LIT}
\end{figure}

We now introduce a key concept for stating our result in the most general form.

\pg{Naive Approach (Optional).}
A naive extension from the static setting leads to a very weak bound. 
\Sps\ $G =G_t = ([N],E)$  for all $t$.
The MSE bound in \cite{jia2023clustered} relies on that for any $i,i',t,t'$ and $r>0$, we have  \begin{align}\label{eqn:111324} {\rm Cov}\lb(\wh \Delta^r_{it},\wh \Delta^r_{i't'}\rb)\ls e^{-{|t-t'|/t_{\rm mix}}}\cdot \mathbbm{1}((i,i') \in E).
\end{align}
This can be {\bf loose} for non-stationary interference graphs.
\Sps\ $i,i'$ were close to each other at time $t_1$ and moved further away from each other since then, and $t_1\ll t'\le t$; see \cref{fig:LIT}.
Intuitively, the covariance between the HT terms should decrease as \( t-t_1\) increases, since the influence of the past diminishes over time.
However, \cref{eqn:111324} only gives a very {\bf loose} bound
\[{\rm Cov} \lb(\wh \Delta^r_{it},\ \wh \Delta^r_{i't'}\rb) = O(1).\]

\pg{Intuition.}
To obtain a tighter bound, we introduce the CIG, and show that the covariance between the HT terms of any individuals $i,i'$ {\bf exponentially} decays in the time elapsed since the last time that $i,i'$ had an edge in the CIG.

To motivate, we first explain how the correlation between the HT terms occurs.
Recall that at the start of each time block $k$, we partition $[N]$ into spatial blocks and assign an arm to each \indep ly.
The state evolution of each $i$ depends on the collection ${\cal E}_{ik}\sse \Pi_k$ of spatial blocks intersected by some $\N_\tau(i)$ where $\tau=t_k,\dots, t_{k+1}-1$.  
For any $i,i'$, if ${\cal E}_{ik} \cap {\cal E}_{i'k}\neq \emptyset$, then both $S_{it}$ and $S_{i't'}$ depend on the arm assigned to a {\bf common} spatial block, which creates a correlation. 

{\new \pg{Clustering-Induced Graphs.}
A vertical design induces a   clustering of space-time into \emph{clusters}, each of which is a product set of a time block and a spatial block.
This clustering induces the following sequence of graphs $(H_k)$. 
For each $k$, we include \( (i, i')\) in the edge set of the $k$-th CIG if there exists a spatial block in \(\Pi_k\) that intersects both \(\mathcal{N}_\tau(i)\) and \( \mathcal{N}_\tau(i')\) for some \( \tau \in [(k-1)\ell, k\ell)\). 
To formalize, recall that \(\mathcal{N}_\tau(i) \subseteq [N]\) denotes the neighbors of individual \( i \) at time \(\tau\).

\bdefn[Clustering-Induced Graph]
Fix a vertical design $(\Pi_k)_{k=1,\dots,T/\ell}$.
The $k$-th {\em clustering-induced graph} (CIG) is denoted \(H_k = ([N], F_k)\), where for $(i,i') \in F_k$ if
\[\exists C \in \Pi_k \ \text{such\ that}\ C\cap \lb( \bigcup_{\tau=t_k}^{t_{k+1}} \mathcal{N}_\tau(i)\rb) \neq \emptyset \ \text{and} \ C \cap \lb( \bigcup_{\tau=t_k}^{t_{k+1}} \mathcal{N}_\tau(i') \rb) \neq \emptyset.\] 
The {\em cluster degree} ${\rm CD}_k(i)$ is the degree of $i$ in $H_k$.
\edefn
We emphasize that if $i=i'$, we trivially have an edge $(i,i)$ in all CIG's.
}

\begin{figure}[h]
\centering
\includegraphics[width=0.9\textwidth]{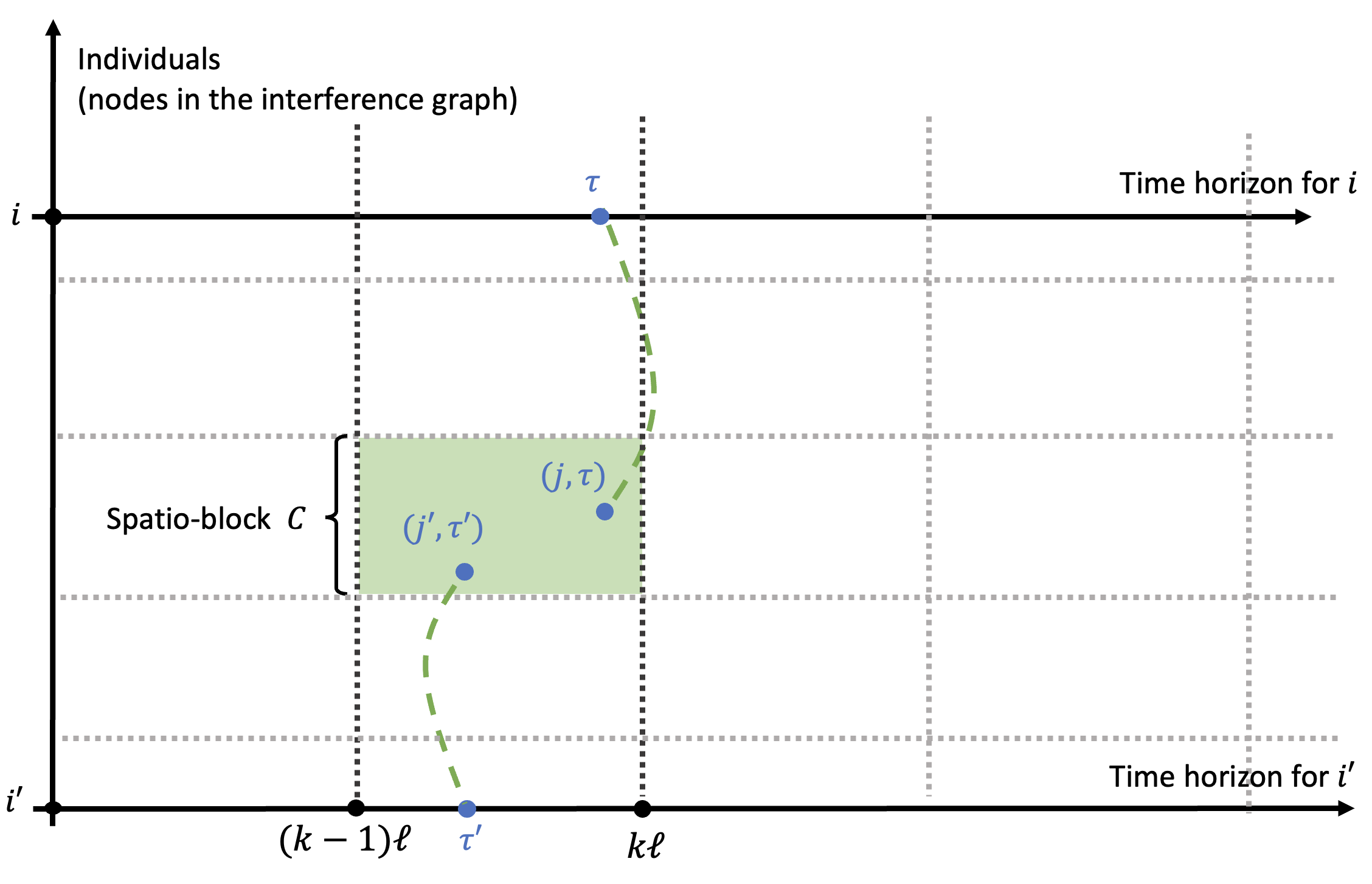}
\caption{{\bf Clustering-induced Graphs.} 
To visualize, we spread out the individuals along a one-dimensional line (the vertical axis). 
Observe from that at time $\tau\in [(k-1)\ell, k\ell]$, there is an edge between $i$ and some individual $j$ in the spatio-block $C\in \Pi_k$, and similarly, $i'$ has an edge with some $j'\in C$ at $\tau'\in [(k-1)\ell, k\ell]$. 
Therefore, we include $(i,i')$ as an edge in the $k$-th CIG.}
\label{fig:LIT}
\end{figure}

\subsection{Illustration: CIG  and Metric-Induced Interference}

For better intuition, we consider the special case where the interference graph is derived from a set of moving individuals in a metric space.  
Given a metric space $({\cal X},d)$ where each point in ${\cal X}$ is called a {\em location}. 
Each individual $i$ has a location $x_t(i)\in {\cal X}$ at each time $t$.
Two individuals can interfere with each other if they are within a given {\rm interference range} $\kappa$. 
Formally, 
\[(i,i') \in E_t \quad {\rm if}\quad  d(x_t(i),\ x_t(i') ) \le \kappa.\]
Equivalently, $G_t$ is the {\em intersection graph} of $\{B(x_t(i), \frac \kappa 2)\}_{i\in [N]}$ where $B(x,\rho)$ denotes the ball with radius $\rho$.
We consider a special case of the vertical design where the space \( \mathcal{X} \) is partitioned into {\em regions} at the start of each time block. 
We randomize at the region level by independently assigning a treatment \( a \in \{0,1\} \) to each region \( R \). 
Specifically, every individual located in region \( R \) at time \( t_k \) receives \( a \) for the {\bf entire duration} of \( [(k-1)\ell, k\ell] \), even if they subsequently move outside the region.
We formalize the above below.

\bdefn[Region-based Vertical Design]
A partition $\Pi$ of $[N]$ is {\em region-based} at time $t$, if there is a partition $\cal R$ of $\cal X$ s.t.
\[\Pi = \{i: x_t(i)\in R\}_{r\in \cal R}.\]
Fix a constant $\ell>0$ for the length of the time blocks. 
A {\em region-based} vertical design is specified by a sequence of partitions $(\Pi_k)_{k=1,\dots,T/\ell}$ of $[N]$, where $\Pi_k$ is a region-based partition at time $t_k$. 
\edefn

\begin{figure}[h]
\centering    \includegraphics[width=0.9\textwidth]
{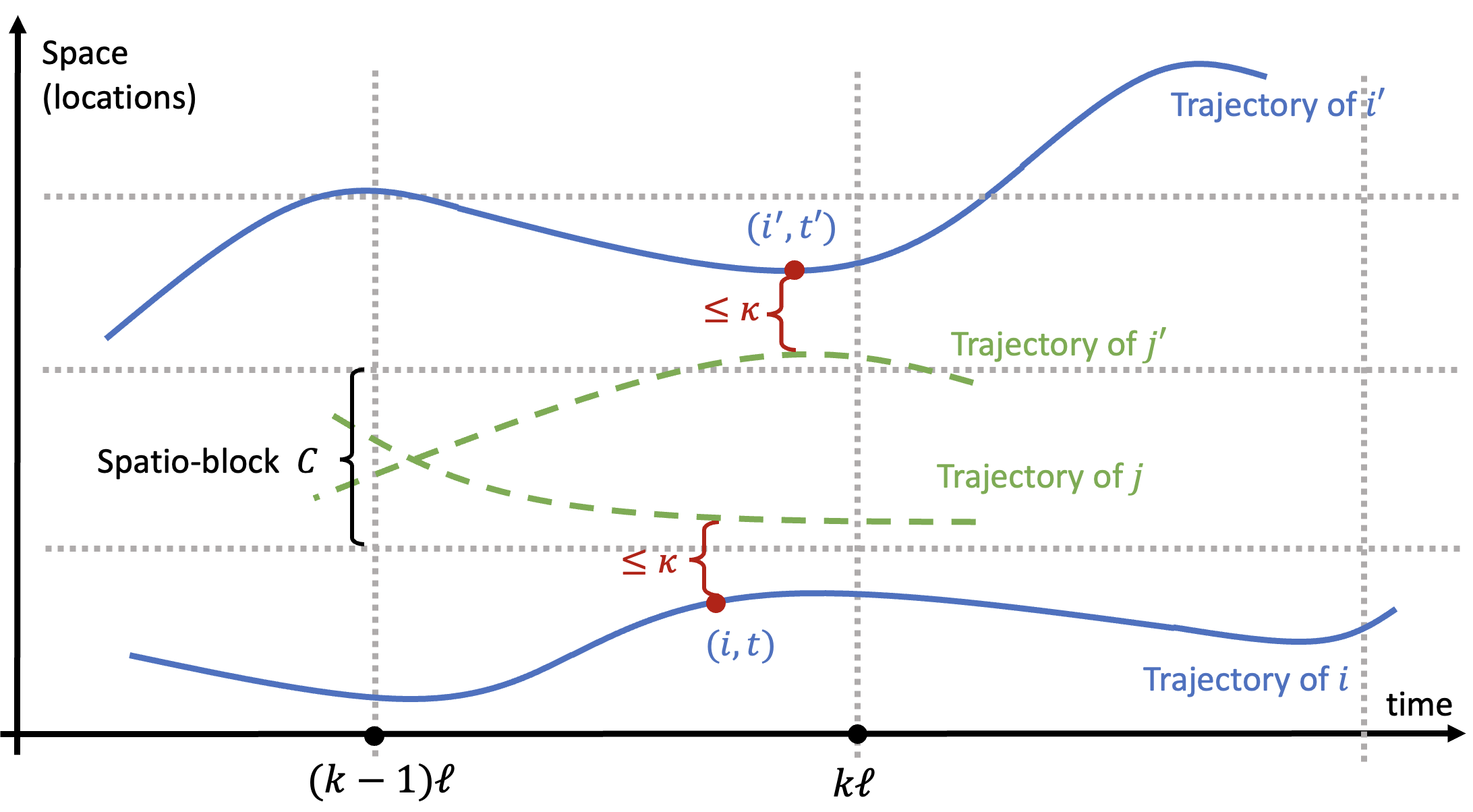}
\caption{{\bf CIG in a  Metric Space.}
To visualize, let us flatten the metric space onto a 1-d line (the vertical axis).
In this figure $(i,i')$ is an edge in the $k$-th CIG. 
In fact, at time $\tau\in [(k-1)\ell, k\ell]$, individual $i$ is within distance $\kappa$ to another individual $j$, who was in spatio-block $C$ at the starting time $(k-1)\ell$ of the time block.
Similarly, $i'$ is within distance $\kappa$ to $j'$, which is also inside $C$ at time $(k-1)\ell$. 
Thus, both $i,i'$ depend on the treatment assigned to $C$ throughout the entire time block, leading to a correlation between their future outcomes and exposure mappings.}
\label{fig:dependence_graph}
\end{figure}

In region-based designs, the notion of $r$-dependence can be characterized geometrically as follows.
For each $i$ and time block $[k\ell, (k+1)\ell]$, the regions that $i$ ``touched'' during $[t-r,t]$ are 
\[{\cal R}_k(it) = \lb\{R\in {\cal R}_k: B \lb(x_t(i), \kappa\rb) \cap R \neq \emptyset {\rm\ for\ some}\ t\in [(k-1)\ell, k\ell) \rb\}.\]


\section{Main Result: A General MSE Bound}
We will first present a covariance bound based on the {\em last interaction time} (LIT), then explain how it leads to an MSE bound.
\subsection{Last Interaction Time (LIT)}
For any fixed \( i, i', t \), their LIT is defined as the endpoint of the most recent time block that (i) begins before \( t \), and (ii) during which \(i,i' \) is an edge in the associated CIG.
\bdefn[Last Interaction Time]   
Fix any $t\in [T]$. 
\Sps\ $t\in [(k-1)\ell,\ k\ell)$ for some $k$.
We define the {\em last interaction time} (LIT) as 
\[\uptau^\star(i,i';t):= k^\star \ell \quad {\rm where} \quad k^\star = \max \{k'\le k: (i,i') \in E_{k'}\}.\]
\edefn


We can now refine the covariance bound, \cref{eqn:111324}, using LIT.

\begin{proposition}[Bounding Covariance Using LIT] 
\label{prop:cov_LIT}
For any $i,i'\in [N]$, $t,t'\in [T]$ \sat ing {\new $t\ge \max\{t',\ \uptau^\star(i,i';t')+r\}$, we have 
\[{\rm Cov}\lb(\wh \Delta^r_{it}, \wh \Delta^r_{i't'}\rb) \ls 
\frac 1{p_{it}^r \cdot p_{i't'}^r} e^{-|t'- \uptau^\star(i,i';t')| / t_{\rm mix}}.\]
}
\end{proposition}




{\new As an important consequence, the above implies that ${\rm Cov}(\wh \Delta_{it}, \wh \Delta_{i't'})$ is only ``non-negligible'' if $t$ and $t'$ lie in the same ``extended'' time block $[(k-1)\ell, k\ell + r]$. 
More precisely: 
\bcoro[Only Close-by Pairs Have Non-negligible Covariance]
\label{coro:only_close_by_pairs_matter}
Fix $i,i'$ and $t'\le t$. 
Denote ${\cal K}_{i,i'} = \{k: (i,i')\in F_k\}$.  
\Sps\ for any $k\in {\cal K}_{i,i'}$, the interval $[(k-1)\ell, k\ell + r]$ contains at most one of $t,t'$. 
Then, 
\[{\rm Cov}\lb(\wh \Delta^r_{it}, \wh \Delta^r_{i't'}\rb) \ls 
\frac 1{p_{it}^r \cdot p_{i't'}^r} e^{-r/ t_{\rm mix}}.\]
\ecoro 
\proof  There are two cases.\\
{\bf Case A:} \Sps\ $t'\notin [k-1)\ell, k\ell + r]$ for any $k\in {\cal K}_{i,i'}$. 
By our assumption that $t\ge t'$, we then have $t\ge \uptau^\star(i,i';,t') + r$. 
Thus, by \cref{prop:cov_LIT}, we have 
\[{\rm Cov}\lb(\wh \Delta^r_{it}, \wh \Delta^r_{i't'}\rb) \le e^{-|\uptau^\star(i,i';,t')|/t_{\rm mix}} \le e^{-r/t_{\rm mix}}.\]
{\bf Case B:} \Sps\ $t'\in [(k-1)\ell, k\ell + r]$ for some $k\in {\cal K}_{i,i'}$.
Then, for any such $k$, we have $t\notin [(k-1)\ell, k\ell + r]$. 
In \parti, $k^*=\max {\cal K}_{i,i'}$, we have $t\notin [(k^*-1)\ell, k^*\ell + r]$.
\OTOH, by the definition of LIT, we have $\uptau^\star(i,i';t')= k^* \ell$.
Thus, the condition 
\[t\ge  \uptau^\star(i,i';t')+r\] 
in \cref{prop:cov_LIT} holds, and we obtain
\[{\rm Cov}\lb(\wh \Delta^r_{it}, \wh \Delta^r_{i't'}\rb) \le e^{-|\uptau^\star(i,i';,t')|/t_{\rm mix}} \le e^{-r/t_{\rm mix}}.\eqno\qed\]
}

\subsection{Main Result}
Recall that the cluster degree  of an individual (w.r.t. time block $k$) is the number of edges incident to it in the $k$-th CIG.
We show that the variance of $\wh \Delta^r$ is essentially $\frac 1{NT}$ times the {\em average cluster degree}, averaged over all $N$ individuals and $T/\ell$ time blocks.

\begin{theorem}[Bounding the Variance Using $\overline{\rm CD}$]
\label{coro:avg_CD}
\Sps\ there is a constant $p^{\rm min}>0$ such that $p_{it}^r\ge p^{\rm min}$ for all $i\in [N],t\in [T]$. 
Then, for any $r,\ell$ we have 
\[{\rm Var}\lb(\wh \Delta^r\rb) \le   \frac {1+\sigma^2}{(p^{\rm min})^2} \cdot \lb(\frac 1{NT} \cdot \frac {(\ell+r)^2}\ell \cdot \overline{\rm CD} + e^{-r/t_{\rm mix}}\rb) \quad {\rm where}\quad \overline{\rm CD} = \frac 1{N \cdot T/\ell}\sum_{i\in [N]} \sum_{k=1}^{T/\ell} {\rm CD}_k(i).\]
\end{theorem}

\proof 
{\new 
We start with decomposing the variance into a sum of covariance terms: 
\begin{align}\label{eqn:072925}
{\rm Var} \lb(\wh \Delta_r\rb) 
&\le \frac 1{N^2 T^2} \sum_{i,i'} \sum_{t'} \sum_{t: t'\le t} {\rm Cov}\lb(\wh \Delta^r_{it}, \wh \Delta^r_{i't'}\rb) \notag\\
&=\frac 1{N^2 T^2} \sum_{i,i'}   M_{i,i'},
\end{align}
where \[M_{i,i'} := \sum_{(t,t'): t'\le t} {\rm Cov}\lb(\wh \Delta^r_{it}, \wh \Delta^r_{i't'}\rb).\] 
We now fix $i,i'$ and bound $M_{i,i'}$ as  
\begin{align*}
M_{i,i'} \le \lb( \sum_{k\in {\cal K}_{i,i'}} \sum_{(t,t')\in [(k-1)\ell,\ k\ell + r]} + 
\sum_{{\rm other}\ t,t' \ {\rm with}\ t'\le t } \rb) {\rm Cov}\lb(\wh \Delta^r_{it}, \wh \Delta^r_{i't'}\rb).
\end{align*}
By  \cref{coro:only_close_by_pairs_matter}, for any $t,t'$ in the second summation, we have 
\[{\rm Cov}\lb(\wh \Delta^r_{it}, \wh \Delta^r_{i't'}\rb) \le \frac 1{p_{it}^r \cdot p_{i't'}^r} \cdot e^{-r/t_{\rm mix}}.\]
Moreover, for each $k\in {\cal K}_{i,i'}$, there are  $O((\ell+r)^2)$ pairs of $t,t'$ in the \corres\ summation, with each covariance at most \[\frac {4(1+\sigma^2)}{p_{it}^r \cdot p_{i't'}^r}\] 
by  \cref{lem:close-by}.
Combining, we obtain
\begin{align*}
M_{i,i'} &\le \sum_{k\in {\cal K}_{i,i'}} (\ell+r)^2 \cdot \frac {4(1+\sigma^2)}{p_{it}^r \cdot p_{i't'}^r} + \sum_{t,t': t'\le t} e^{-r/t_{\rm mix}}\\
& \le \sum_{k\in {\cal K}_{i,i'}} (\ell+r)^2 \cdot \frac {4(1+\sigma^2)}{p_{it}^r \cdot p_{i't'}^r} + \sum_{t,t': t'\le t} \frac 1{p_{it}^r \cdot p_{i't'}^r} e^{-r/t_{\rm mix}}\\
&\le \frac {4(1+\sigma^2)}{(p^{\rm min})^2} \lb((\ell+r)^2 \cdot |{\cal K}_{i,i'}| + T^2 e^{-r/t_{\rm mix}} \rb).
\end{align*}
Substituting the above into \cref{eqn:072925} and noting that 
\[2\sum_{i,i'} |{\cal K}_{i,i'}| = \sum_{i,i'} {\rm CD}_k(i) =  N\cdot \frac T\ell \cdot \overline {\rm CD},\]
we obtain
\begin{align*}
{\rm Var} \lb(\wh \Delta_r\rb) 
&\le \frac {4(1+\sigma^2)}{(p^{\rm min})^2} \cdot \frac 1{N^2 T^2}  \lb( \sum_{i,i'\in [N]}\lb((\ell+r)^2 \cdot |{\cal K}_{i,i'}| + T^2 e^{-r/t_{\rm mix}}\rb) \rb)\\
&\le \frac {2(1+\sigma^2)}{(p^{\rm min})^2} \cdot \frac 1{N^2 T^2} \lb( (\ell+r)^2 \cdot N\cdot \frac T\ell \cdot \overline{\rm CD} + N^2 T^2 e^{-r/t_{\rm mix}}\rb) \\
&\le \frac {2(1+\sigma^2)}{(p^{\rm min})^2} \cdot \lb(\frac 1{NT} \cdot \frac {(\ell+r)^2}\ell \cdot \overline{\rm CD} + e^{-r/t_{\rm mix}}\rb).\tag*\qed
\end{align*}

To further simplify, we choose \( r, \ell\) to be logarithmic in $NT$, so that the exponential term becomes lower-order: 
\bcoro
With $\ell=r=2t_{\rm mix} \log NT$, we have 
\[{\rm Var} \lb(\wh \Delta_r\rb) \ls \frac 1{(p^{\rm min})^2} \cdot \frac {t_{\rm mix} \log NT}{NT} \cdot \overline{\rm CD}.\]
\ecoro 
We will soon explain how to further simplify and interpret the above.
}

\subsection{Minimal Exposure \Prb}
To further simplify the above, we bound the exposure \prb\ in terms of the cluster degrees.
\bprop[Lower Bound for the Exposure \Prb] 
For any $i,t,r$, we have 
\[p_{it}^r \ge 2^{-\sum_{k: [t_{k-1},t_k] \cap [t-r, t]\neq \emptyset} {\rm CD}_k(i)}.\]
\eprop
To see this, observe that the right-hand side corresponds to the \textbf{exact} exposure probability $\prb$ in the most extreme case -i.e., when every neighbor of $i$ in the graph $H_k$ lies in a distinct spatial block, for every $k$ with $B_k \cap [t - r, t] \neq \emptyset$.
(Of course, in many practically motivated interference structures, the actual exposure probability \prb\ can be much higher, as we will see in the next section.)

\brmk This bound is worst-case - the identify holds only when every neighbor lies in a different spatial block.
\ermk 

As an important case, if $\ell=r$, there are at most two time blocks that intersect $[t-r,t]$, so 
\[p_{it}^r \ge 2^{-2{\rm CD}_{\rm max}}\quad {\rm where}\quad {\rm CD}_{\rm max} = \max_{i\in [N], 1\le k\le T/\ell} {\rm CD}_k(i).\] 
We obtain the following by observing that the bias bound \cref{prop:bias} becomes a lower-order term when $r$ is sufficiently large.

\begin{theorem}[Bounding MSE Using ACD] 
\label{thm:mse}
\Sps\ $\ell = r =  2t_{\rm mix}\log NT$. Then,
\[{\rm MSE} = O\lb(\frac {t_{\rm mix} \log (NT)}{NT} \cdot 2^{4 {\rm CD}_{\rm max}} \cdot \overline{\rm CD}\rb).\]
\end{theorem}

\brmk As a sanity check, consider the \emph{no-interference} setting, where each individual forms its own cluster.
Since the interference graph has no edges, the corresponding cluster graph also contains no edges apart from self-loops.  
Therefore, we have ${\rm CD}_k(i) = 1$ for all $i$ and $k$, and 
\cref{thm:mse} gives a variance bound of $\tilde O(\frac 1{NT})$.

When the graph changes, the CD's can increase, as they are defined based on the union of connections across all rounds within a time block.
Furthermore, the faster the graph changes, the larger the resulting CD's tend to be.
\qed 
\ermk 

\su{ADD A SECTION: Discuss implications of the main theorem in the metric  setting}

\su{Special case: two extremes  (1) draw a new Erdos Renyi random graph (2) static. 
Interpolate: each edge die/born with some \prb\ $p$. Add some ``analytical tools'' to analyze this dependence $p$ }

\su{Also: Beef up practical insights regarding spatial non-uniformity. 
How should we change the spatio-partitioning, by accounting for how interference graphs changes (peak/off peak for example), and both $t_{mix}$. 
}

\bibliography{ref}
\bibliographystyle{plainnat}

\appendix
\section{Bias analysis}
\label{sec:bias}
For any event $\mathcal{F}_ W$-measurable event $A$, denote by $\mathbb P_A$, $\mathbb E_A$,  and $\mathrm{Cov}_A$ probability, expectation, and covariance conditioned on $A$.
If $A=\{w\}$ is singleton, we use the subscript $w$.
We first show that if $w,w'$ are  identical in a neighborhood of $it$, then the conditional \distr s of $Y_{it}$ on $w,w'$ are also close. 

\begin{lemma}[Decaying Temporal Interference]
\label{prop:TVdecay}
Consider any $t,r\in [T]$ with $1\leq r<t\leq T$ and $i\in [N]$.
Suppose $w,w'\in \{0,1\}^{N\times T}$ are identical on $\N^r(it)$.
Then, \[d_{\rm TV} \lb(\mathbb P_w[Y_{it}\in\cdot], \ \mathbb P_{w'}[Y_{it}\in\cdot]\rb) \le e^{-r/t_{\rm mix}}.\]
\end{lemma}
\proof For any $i,t$, we denote $f_{it}= \mathbb P_w[S_{it} \in\cdot]$ and $f'_{it}= \mathbb P_{w'}[S_{it}\in\cdot]$.
Then, by the Chapman–Kolmogorov equation, for any $\tau$ we have 
\[f_{i,\tau+1} = f_{i\tau} P^{w_{\N_\tau(i)},\tau}_{i\tau}\quad \text{and} \quad f'_{i, \tau+1} = f'_{i\tau} P^{w'_{\N_\tau(i)},\tau}_{i\tau}.\]
Thus, if $t-r\le \tau\le t$, 
\begin{align*}
d_{\rm TV}(f_{i, \tau+1},\ f'_{i, \tau+1})&= d_{\rm TV}\lb(f_{i\tau} P^{w_{\N_\tau(i)},\tau}_{i\tau},\ f'_{i\tau} P^{w'_{\N_\tau(i)},\tau}_{i\tau}\rb)\\
&= d_{\rm TV}\lb(f_{i\tau} P^{w_{\N_\tau(i)},\tau}_{i\tau},\ f'_{i\tau} P^{w_{\N_\tau(i)},\tau}_{i\tau}\rb)\\
&\le e^{-1/t_{\rm mix}} \cdot d_{\rm TV}\lb(f_{i\tau}, f'_{i\tau}\rb),
\end{align*}
where we used $w_{\N^r(it)} = w'_{\N^r(it)}$ in the second equality and the rapid mixing \assu\ in the inequality.
Applying the above for all $\tau =t-r,\dots,t$, we conclude that
\begin{align*}
d_{\rm TV} \lb(\mathbb P_w[Y_{it}\in\cdot],\ \mathbb P_{w'}[Y_{it}\in\cdot]\rb)  
&\le d_{\rm TV} \lb(f_{it},f'_{it}\rb)\\ &\le e^{-r/ t_{\rm mix}} \cdot d_{\rm TV}(f_{i,t-r},\ f'_{i,t-r})\\
&\le e^{-r/ t_{\rm mix}},
\end{align*}
where the first inequality is because $Y_{it}= \mu_{it}(S_{it}, w_{\N_t(i), t}) + \eps_{it}$ where  $\eps_{it}$ is an external error and $|\mu_{it}(\cdot,\cdot)| \le 1$.
\qed 

We next use the above to show that  the bias is $e^{-r/t_{\rm mix}}$. 

\noindent{\bf Proof of \cref{prop:bias}.} 
Fix any $a\in \{0,1\}$, $r>0$ and $(i,t)\in [N]\times [T]$.
Consider any $w$ with $X_{ita}^r(w)=1$, or equivalently, $w_{\N^r(it)}=a\cdot{\bf 1}$.
So, by \cref{prop:TVdecay} (with $a\cdot {\bf 1}$ in the role of $W'$), we have
\begin{align}\label{eqn:053125}
d_{\rm TV} \lb(\mathbb P[Y_{it}\in\cdot\mid W=w],\ \mathbb P[Y_{it}\in\cdot\mid W= {a\cdot \bf 1}]\rb) \le e^{-r/t_{\rm mix}}.
\end{align}
Observing that $\mathbb P[Y_{it}\in\cdot\mid X_{ita}^r(W)=1]$ is a convex combination of $\mathbb P_W[Y_{it}\in\cdot\ ]$ terms over all $W$'s with $X_{ita}^r(W)=1$, we have
\[d_{\rm TV} \lb(\mathbb P[Y_{it}\in\cdot\mid X_{ita}^r(W)=1],\ \mathbb P[Y_{it}\in\cdot\mid W= {a\cdot \bf 1}]\rb) \le e^{-r/t_{\rm mix}}.\]
Now we use \cref{eqn:053125} to bound the bias:
Conditioning on $X_{ita}^r$, we have
\begin{align*}
\ho{E}\lb[\wh Y_{ita}^r\rb] 
&= \ho{E}\lb[\frac{X_{ita}^r}{p_{it}^r} Y_{it} \Bigm| X_{ita}^r = 1\rb] \ho{P}\lb[X_{ita}^r = 1\rb] + \ho{E}\lb[\frac{X_{ita}^r}{p_{it}^r} Y_{it} \Bigm| X_{ita}^r = 0\rb] \ho{P}\lb[X_{ita}^r = 0\rb] \\
&= \ho{E}\lb[\frac{X_{ita}^r}{p_{it}^r} Y_{it}\Bigm| X_{ita}^r = 1\rb] p_{it}^r +0 \\
& = \ho{E} \lb[Y_{it} \Bigm| X_{ita}^r =1\rb].
\end{align*}
Combining the above with \cref{eqn:053125}, we have
\begin{align*}
&\quad \lb|\ho{E}\lb[\wh Y_{ita}^r\rb] - \ho{E}\lb[Y_{it}\mid W =a\cdot {\bf 1}\rb]\rb| \\
&= \lb|\ho{E} \lb[ Y_{it} \mid X_{ita}^r(W)=1\rb] - \ho{E}\lb[Y_{it}\mid W = a\cdot {\bf 1}\rb] \rb| \\
&\le e^{-r/t_{\rm mix}}.
\end{align*}
Finally, by triangle \ineq, we conclude that
\begin{align*}
\lb|\ho{E} \lb[\wh\Delta^r\rb] - \Delta\rb| 
& \le \frac 1{NT}\sum_{(i,t)\in [N]\times [T]} \lb|\Delta_{it} - \ho{E}\lb[\wh \Delta_{it}^r\rb]\rb| \\
& \le \frac 1{NT} \sum_{a\in \{0,1\}}\sum_{i,t}\lb|\ho{E}\lb[\wh Y_{ita}^r\rb] - \ho{E}\lb[Y_{it}\mid W =a\cdot {\bf 1}\rb] \rb|\\
&\le 2e^{-r/t_{\rm mix}}.\tag*\qed
\end{align*}

\section{Proof the Covariance Bound (\cref{prop:cov_LIT})} \label{sec:var}
We start with a bound that holds for all pairs $it, i't'$ of space-time points.
It follows from the Cauchy-Schwarz \ineq\ and the boundedness of the \rv s.

\begin{lemma}[Covariance of HT Terms]\label{lem:close-by}  
For any $r\ge 0$, $i,i'\in [N]$,  $t,t'\in [T]$ and design $\cal D$, we have
\[{\rm Cov}\lb(\wh\Delta_{it}^r,\ \wh\Delta_{i't'}^r\rb) \le  \frac{4(1+\sigma^2)}{p_{it}^r\cdot p_{i't'}^r}.\]
\end{lemma}
\proof Expanding the definition of $\wh \Delta_{it}^r$, we have  
\begin{align}\label{eqn:121723}
{\rm Cov} \lb(\wh\Delta_{it}^r,\ \wh\Delta_{i't'}^r\rb) 
& = {\rm Cov}\lb(\frac{X_{it1}^r}{p_{it}^r} Y_{it} - \frac{X_{it0}^r}{p_{it}^r} Y_{it},\  \frac{X_{i't'1}^r}{p_{i't'}^r}Y_{i't'}- \frac{X_{ i't'0}^r}{p_{i't'}^r}Y_{i't'}\rb)\notag \\
&\le \sum_{a,a'\in \{0,1\}} \lb|{\rm Cov} \lb(\frac{X_{ita}^r}{p_{it}^r} Y_{it},\ \frac{X_{i't'a'}^r}{p_{i't'}^r} Y_{ i't'}\rb)\rb| \notag\\
&\le \frac 1{p_{it}^r}\frac 1{p_{i't'}^r} \sum_{a,a'\in \{0,1\}}  \lb|{\rm Cov} \lb(X_{ita}^r Y_{it},\ X_{i't'a'}^r Y_{ i't'}\rb)\rb| \notag\\
&\le \frac 1{p_{it}^r}\frac 1{p_{i't'}^r} \sum_{a,a'\in \{0,1\}} \sqrt {\ho{E}\lb[(X_{ita}^r Y_{it})^2\rb]} \cdot \sqrt{\ho{E}[(X_{i't'a'}^r Y_{i't'})^2]} \notag\\
&\le \frac {4(1+\sigma^2)}{p_{it}^r\cdot p_{i't'}^r},
\end{align}
where the second last inequality follows from the Cauchy-Schwarz inequality and last \ineq\ follows since $\ho{E}[(X_{ita}^r Y_{it})^2]\le 1+\sigma^2$. \qed 

The above bound alone is not sufficient for our analysis, as it does not take advantage of the rapid mixing property.
The rest of this section is dedicated to showing that  HT terms decays {\bf exponentially} in the temporal distance.
To this end, we will decompose the covariance as follows.

\bprop[Law of Total Covariance]
For any \rv s $U,V,W$, we have \begin{align*}
{\rm Cov}(U,V) = \ho{E}[{\rm Cov}(U,V\mid W)] + {\rm Cov}(\ho{E}[U\mid W],\ \ho{E}[V\mid W]).
\end{align*}
\eprop

We will apply the above as follows. 
By the Law of Total Covariance (LOTC), for any individuals $i,i'\in [N]$, rounds $t,t'\in [T]$ and radius $r\ge 0$ for the HT estimator, we have
\begin{align}\label{eqn:LOTC} 
{\rm Cov}\lb(\wh\Delta_{it}^r, \wh\Delta_{i't'}^r\rb) = \ho{E}\lb[{\rm Cov} \lb(\wh \Delta_{it}^r, \wh \Delta_{i't'}^r\middle|\ W\rb)\rb] + {\rm Cov}\lb(\ho{E}\lb[\wh \Delta_{it}^r\middle|\ W\rb],\ \ho{E}\lb[\wh \Delta_{i't'}^r\middle|\ W\rb]\rb).
\end{align}

We will bound these two terms separately in the next two subsections. 
Throughout, we will sometimes write $\ho{P}_w = \ho{P}[\cdot \mid W=w]$ and ${\rm Cov}_w (\cdot,\cdot) = {\rm Cov}(\cdot,\cdot \mid W=w)$.

\subsection{Bounding the First Term in LOTC}
We first show that if the realization of one random variable has little impact on the (conditional) distribution of another random variable, then they have low covariance.

\begin{lemma}[Low Interference in Conditional Distribution Implies Low Covariance]\label{lem:low_cov}
Let $U,V$ be two random variables and $g,h$ be real-valued functions defined on their respective realization spaces.
If for some $\delta>0$,  
\[d_{\rm TV}(\mathbb P[U\in\cdot\mid V],\ \mathbb P[U\in\cdot])\le \delta \quad V\text{-almost surely },\]
then, 
\[{\rm Cov}(g(U),h(V)) \le \delta \cdot \|h(V)\|_1 \cdot  \|g(U)\|_\infty.\]
\end{lemma}
\proof Denote by $\mu_{U,V},\mu_U,\mu_V,\mu_{U\mid V=v}$ the probability measures of $(U,V)$, $U$, $V$, and $U$ conditioned on $V=v$, respectively. We then have
\begin{align*} 
|{\rm Cov}(g(U),h(V))| &= |\ho{E}\lb [g(U) h(V)\rb] - \ho{E}[g(U)] \ho{E} [h(V)]| \notag\\ 
& =\lb|\int_v h(v) \lb(\int_u g(u) \lb(\mu_{U\mid V=v}(du) - \mu_U(du) \rb)\rb)\mu_V(dv)\rb|\\
&\le \int_v |h(v)|  \cdot \|g(U)\|_\infty \cdot 
d_{\rm TV}(\mathbb P(U\in\cdot\mid V=v), \mathbb P(U\in\cdot))\ \mu_V(dv)\\
&\le \|h(V)\|_1\cdot  \|g(U)\|_\infty \cdot  \delta. \tag*\qed
\end{align*}

To use the above, we will need to bound the TV-distance between the conditional \distr s.
The key is to analyze how much $\ho{P}[S_{it}\in \cdot\mid S_{i't'}=s]$ is affected by $s$. 
An obvious bound is \[d_{\rm TV}(\ho{P}[S_{it}\in \cdot\mid S_{i't'}=s],\ \ho{P}[S_{it}\in \cdot\mid S_{i't'}=s'])\le e^{-|t-t'|/t_{\rm mix}}.\]
However, this is still quite loose. 
For example, for $t=t'$, this bound
which is {\bf hard to analyze} directly since $S_{it},S_{i't'}$. 
To mitigate this, we {\bf backtrack} to the last time when $i$ and $i'$ interact, and argue that the interference decays  since then exponentially. 
More precisely, this involves writing the target conditional \prb\ as a \cvx\ combination of terms of the form $\ho{P}[S_{it}\in \cdot\mid S_{i,\uptau^\star}\in \cdot]$, and use the rapid mixing \cond\ (\cref{assu:rapid_mixing}) to bound each of them. 
To formalize, we need:

\begin{lemma}[Decomposition of conditional \prb]
\label{lem:cond_ind}  
Let $Z$ be a \rv\ taking values on a countable set $\cal Z$ and $A,B$ be $Z$-\msb\ events. 
\Sps\ $A\ind B \mid Z.$ Then, 
\[\ho{P}[A\mid B] = \sum_{z\in \cal Z} \ho{P}[A\mid Z=z] \cdot \ho{P}[Z=z\mid B].\]
\end{lemma}
\proof By the definition of conditional \prb, we have
\begin{align}\label{eqn:072625}
\sum_{z\in \cal Z}\ho{P}[A\mid Z=z]\cdot \ho{P}[Z=z\mid B]
&= \sum_{z\in \cal Z} \ho{P}[A\mid Z=z] \cdot \frac{\ho{P}[\{Z=z\}\cap B]}{\ho{P}[Z=z]} \cdot \frac{\ho{P}[Z=z]}{\ho{P}[B]} \notag \\
&=  \sum_{z\in \cal Z} \ho{P}[A\mid Z=z] \cdot \ho{P}[B\mid Z=z] \cdot \frac{\ho{P}[Z=z]}{\ho{P}[B]}
\end{align}
Since $A\ind B\mid Z$, we have  
\begin{align*}
\eqref{eqn:072625} & =\frac{\sum_{z\in \cal Z} \ho{P}[A\cap B \mid Z=z] \cdot \ho{P}[Z=z]}{\ho{P}[B]}\\
&= \frac{\sum_{z\in \cal Z} \ho{P}[A\cap B] \cdot \ho{P}[Z=z]}{\ho{P}[B]} \\
&=\ho{P}[A\mid B]. \tag*\qed 
\end{align*}

We next show that the state distribution of an individual at time \( t \) is largely insensitive to the state of another individual at time \( t' \), provided their last interaction occurred sufficiently long ago. 
This follows from the result above by treating \( S_{i,\uptau^\star} \) as \( Z \). 
For convenience, write $\uptau^\star = \uptau^\star(i, i'; t)$.

\begin{lemma}[Exponential Decay in the TV-Distance]
\label{lem:cov_states}
\Sps\ $i,i'\in [N]$ and {\new $t\ge \max\{t',\uptau^\star(i,i';t')\}$.}
Let $B$ be an $S_{i't'}$-\msb\ event. 
Then, for any $W\in \{0,1\}^{N\times T}$, we have
\begin{align}\label{eqn:082824} 
d_{\rm TV}(\ho{P}_w[S_{it}\in \cdot \mid S_{i't'}\in B],\ \ho{P}_w[S_{it}\in \cdot\ ])\le e^{-|t-\uptau^\star(i,i';t')|/ t_{\rm mix}}. 
\end{align}
\end{lemma} 
\proof By the definition of $\uptau^\star:=\uptau^\star(i,i';t')$, we crucially have 
\[{\blue S_{it} \ind S_{i't'} \mid \{S_{i, \uptau^\star},S_{i', \uptau^\star}\}}.\] 
So by \cref{lem:cond_ind},  
\begin{align}\label{eqn:082924}
&\quad\ \ho{P}[S_{it}\in \cdot \mid S_{i't'}\in B]\notag \\
&= \sum_{s,s'\in \cal S} \ho{P}[S_{it} \in \cdot \mid (S_{i,\uptau^\star}, S_{i',\uptau^\star})=(s,s')]\cdot \ho{P}[(S_{i,\uptau^\star}, S_{i',\uptau^\star})=(s,s')\mid S_{i't'} \in B].
\end{align}
Moreover, for any $s,s'\in\cal S$, the rapid mixing \assu\ gives 
\begin{align*}
&\quad\ d_{\rm TV} \lb(\ho{P}[S_{it} \in \cdot \mid S_{i,\uptau^\star}=s,\ S_{i',\uptau^\star}=s'],\ \ho{P}[S_{it}\in \cdot\ ]\rb) \\
&= d_{\rm TV} \lb(\ho{P}[S_{it} \in \cdot \mid S_{i,\uptau^\star}=s],\ \ho{P}[S_{it}\in \cdot\ ]\rb)\\
&\le e^{-|t-\uptau^\star| / t_{\rm mix}}.
\end{align*}
Therefore, 
\begin{align*}
\eqref{eqn:082924} &\le \sum_{s,s'\in \cal S} e^{-|t-\uptau^\star| / t_{\rm mix}} \cdot \ho{P}[(S_{i,\uptau^\star}, S_{i',\uptau^\star})=(s,s')\mid S_{i't'} \in B]\\
&=e^{-|t-\uptau^\star| / t_{\rm mix}}.\tag*\qed
\end{align*}

\begin{lemma}[Covariance in the Outcomes] \label{lem:cov_outcomes}
For any $w\in \{0,1\}^{N\times T}$, $i,i'\in [N]$ and $t' \le t$, we have 
\[{\rm Cov}(Y_{it}, Y_{i't'}\mid W=w) \le e^{-{|t-\uptau^\star(i,i';t)| /{t_{{\rm mix}}}}}.\]
\end{lemma} 
\proof 
Recall that $\epsilon_{it} =Y_{it}-\mu_{it}(S_{it},\ W_{\N_t(i),t})$, so
\begin{align*}
{\rm Cov}\lb(Y_{it},\ Y_{i't'}\mid W=w\rb)
&= {\rm Cov}(\mu_{it}(S_{it},W_{\N_t(i),t}),\ \mu_{i't'}(S_{i't'},W_{\N_t(i'),t'})\mid W=w)\\
& \quad + {\rm Cov}(\eps_{i't'},\ \mu_{it}(S_{it},W_{\N_t(i),t})\mid W=w) \\
&\quad +{\rm Cov}(\mu_{i't'}(S_{i't'},W_{\N_t(i'),t'}),\ \eps_{it}\mid W=w) \\
& \quad + {\rm Cov} (\epsilon_{i't'},\ \epsilon_{it}\mid W=w).
\end{align*}
The latter three terms are $0$ by the exogenous noise assumption.
To bound the first term, note that  by \cref{lem:cov_states}, for any $S_{i't'}$-\msb\ event $B$, we have
\begin{align}\label{eqn:121423b} 
&\ \quad d_{\rm TV}\lb(\ho{P}_w \lb[(S_{it}, W_{\N_t(i), t})\in\cdot\ \rb],\ \ho{P}_w \lb[(S_{it}, W_{\N_t(i),t})\in\cdot \mid S_{i't'}\rb]\rb)\notag\\
& = d_{\rm TV}\lb(\ho{P}_w \lb[S_{it}\in\cdot\ \rb],\ \ho{P}_w \lb [S_{it}\in\cdot\mid S_{i't'}\rb]\rb) \notag\\
&\le e^{-{|t-\uptau^\star(i,i';t)| /{t_{{\rm mix}}}}}.
\end{align}
To conclude, apply \cref{lem:low_cov} with $(S_{ i't'},W_{ i't'})$ in the role of $V$, with $(S_{it},W_{it})$ in the role of $U$, with $\mu_{it},\mu_{ i't'}$ in the role of $g,h$, and with $e^{-|t-\uptau^\star|/t_{\rm mix}}$ in the role of $\delta$, and noting that $\|g\|_\infty, \|h\|_\infty\le 1$. 
\cref{eqn:121423b} then immediately leads to the claimed \ineq.
\hfill\qed




As the final building block, we  simplify the covariance by conditioning on  the exposure mappings.

\begin{lemma}[Expectation of  Conditional Covariance]
\label{lem:exp_cond_cov} 
For any $a,a'\in \{0,1\}$, write 
\[E_{a,a'}= \lb\{w\in \{0,1\}^{N\times T}: X_{ita}^r(w)=X_{i't'a'}^r(w)=1\rb\}.\]
Then, for any $i,i'\in [N]$ and $t'\le t$ such that $t\ge \uptau^\star(i,i';t') +r$, we have
\[\ho{E}\lb[{\rm Cov} \lb(\wh \Delta_{it}^r, \wh \Delta_{i't'}^r\middle |\ W\rb)\rb] 
= \sum_{(a,a')\in \{0,1\}^2} (-1)^{a+a'} \cdot \ho{E}\lb[{\rm Cov}\lb(Y_{it}, Y_{i't'}\mid W\rb)\mid W\in E_{a,a'}\rb].\]
\end{lemma}
\proof Fix any $a,a'\in \{0,1\}$. 
Since $p_{it}^r, p_{i't'}^r$ are constants, we have
\begin{align}\label{eqn:052425c}
\ho{E}\lb[{\rm Cov}\lb(\frac{X_{ita}^r}{p_{it}^r} Y_{it}, \frac{X_{i't'a'}^r}{p_{i't'}^r} Y_{i't'}\ \middle |\ W\rb)\rb] 
= \frac 1{p_{it}^r \cdot p_{i't'}^r} \ho{E}\lb[{\rm Cov} (X_{ita}^r Y_{it},\ X_{i't'a'}^r Y_{i't'}\mid W)\rb].
\end{align}
Note that $X_{ita}^r,X_{i't'a'}^r$ are functions of $W$ (and hence deterministic once $W$ is fixed), so
\begin{align}\label{eqn:052425d}
\eqref{eqn:052425c} 
&= \frac 1{p_{it}^r p_{i't'}^r} \ho{E}\lb[X_{ita}^r(W)\cdot X_{i't'a'}^r(W) \cdot {\rm Cov} (Y_{it}, Y_{i't'}\mid W)\rb].
\end{align}
Since $X_{ita}^r, X_{i't'a'}^r\in \{0,1\}$, the term ${\rm Cov} (X_{ita}^rY_{it}, X_{i't'a'}^r Y_{i't'}\mid W)$ is non-zero {\bf only} when $W\in E_{a,a'}$.
Moreover, by the definition of $\uptau^\star$ and our assumption that $t\ge \uptau^\star +r$, we have 
\[X_{ita}^r\ind X_{i't'a'}^r.\] 
Thus,
\[\ho{P}[W\in E_{a,a'}] = \ho{P}[X_{ita}^r=1]\cdot \ho{P}[X_{i't'a'}^r=1] = p_{it}^r \cdot p_{i't'}^r,\]
and
\begin{align}\label{eqn:052425e} 
\eqref{eqn:052425d}
&= \frac 1{p_{it}^r \cdot p_{i't'}^r} \cdot \ho{P}[X_{ita}^r(W) = X_{i't'a'}^r(W) = 1] \cdot \ho{E}[{\rm Cov} (Y_{it},Y_{i't'}\mid W)\mid X_{ita}^r = X_{i't'a'}^r = 1]\notag \\
& = \ho{E}\lb[{\rm Cov}\lb(Y_{it}, Y_{i't'}\mid W\rb)\mid W\in E_{a,a'}\rb].
\end{align}
Applying the above for all $a,a'$ and using the definition of $\wh \Delta_{it}^r$, we conclude that 
\[\ho{E}\lb[{\rm Cov} \lb(\wh \Delta_{it}^r, \wh \Delta_{i't'}^r\middle|\ W\rb)\rb]  
= \sum_{(a,a')\in \{0,1\}^2} (-1)^{a+a'} \cdot \ho{E}\lb[{\rm Cov}\lb(Y_{it}, Y_{i't'}\ \middle |\ W\rb) \mid W\in E_{a,a'}\rb].\eqno\qed
\]

Now we are able to bound the first term in \cref{eqn:LOTC} (the decomposition of covariance based on the LOTC).

\begin{proposition}[Bounding the Expected Conditional Covariance]
\label{prop:red_dot} 
For any $i,i'\in [N], 1\le t\le t'\le T$ such that {\new $t\le \uptau^\star(i,i';t) + r$,} we have  
\[\ho{E}\lb[{\rm Cov}\lb(\wh \Delta_{it}^r, \wh \Delta_{i't'}^r\ \middle |\ W\rb)\rb] 
{\new \le e^{-|t-\uptau^\star(i,i';t')|/ t_{\rm mix}} }.\]
\end{proposition}
\proof Fix $a,a'\in \{0,1\}$ and $w\in \{0,1\}^{N\times T}$. 
Observe that if $X_{ita}^r(w) = 0$ or $X_{i't'a'}^r (w) = 0$, then 
\[{\rm Cov} \lb(X_{ita}^r Y_{it},\ X_{i't'a'}^r Y_{i't'}\middle |\ W=w\rb) = 0.\] 
Now \sps\ $X_{ita}^r(w) = X_{i't'a'}^r(w) = 1$. 
In this case, we have 
\[{\rm Cov} \lb(X_{ita}^r Y_{it},\ X_{i't'a'}^r 
 Y_{i't'}\mid W=w\rb) = {\rm Cov}(Y_{it},\ Y_{i't'}\mid W=w).\]
We now use \cref{lem:low_cov} to bound the above.
By \cref{lem:cov_outcomes}, for any $Y_{i't'}$-\msb\ events $A,B$, we have
\[d_{\rm TV}\lb(\ho{P}_w[Y_{it}\in \cdot\mid Y_{i't'}\in A],\ \ho{P}_w[Y_{it}\in \cdot\mid Y_{i't'}\in B]\rb)\le {\new  e^{-|t-\uptau^\star(i,i';t')|/t_{\rm mix}}}.\]
Thus, by \cref{lem:low_cov}, for any $w$ we have
\[\lb|{\rm Cov}(Y_{it}, Y_{i't'}\mid W=w)\rb| \le {\new e^{-|t-\uptau^\star(i,i';t')|/t_{\rm mix}}}.\] 
Thus,
\[\lb| \ho{E} [{\rm Cov}(Y_{it}, Y_{i't'}\mid W) \mid W\in E_{a,a'}] \rb| \le {\new e^{-|t-\uptau^\star(i,i';t')|/t_{\rm mix}}}.\]
Finally, by  \cref{lem:exp_cond_cov} and {\new our assumption that $r\le \uptau^\star(i,i';t)$}, we conclude that 
\[\ho{E}\lb[{\rm Cov} \lb(\wh \Delta_{it}^r, \wh \Delta_{it'}^r\middle |\ W\rb)\rb] 
\le \sum_{(a,a')\in \{0,1\}^2}  \lb|\ho{E}\lb[{\rm Cov}\lb(Y_{it}, Y_{it'}\mid W\rb)\mid W\in E_{a,a'}\rb]\rb|
\le {\new 4 e^{-|t-\uptau^\star(i,i';t')|/t_{\rm mix}}}.\eqno\qed\]

\subsection{Bounding the Second Term in LOTC}
The following says that if $it$ and $i't'$ are far apart - in that there is a pair of large neighborhoods that intersect no common clusters, then their expected HT terms conditional on the treatment assignment $W$ have low covariance.

\su{Add another prop to consider  the other case, where $t'$ in in $\tau^*$'s time block, but $t-\tau^*$ is large. 
Proof idea:  only  decomp the cond exp of $t$. Then use the fact that $W_{recent}\ind W_{before\ \tau^*}$. Then CS ineq.
}
\begin{proposition}
[Covariance of Conditional Expectation]
Fix $i,i'\in [N]$, $t'\le t$ and $r\ge 0$ with $t'-r \ge \uptau^\star:=\uptau^\star(i,i';t')$. 
Then,  
\[{\rm Cov}\lb(\ho{E}\lb[\wh \Delta_{it}^r\middle |\ W\rb],\ \ho{E} \lb[\wh \Delta_{i't'}^r\middle |\  W\rb]\rb)\le \frac {12 e^{-|t'-\uptau^\star| / t_{\rm mix}}}{p_{it}^r\cdot p_{i't'}^r}.\] 
\end{proposition}

\proof  
Fix $a,a'\in \{0,1\}$. 
Write 
\begin{align}\label{eqn:052425f}
&\quad \ {\rm Cov}\lb(\ho{E} \lb[\frac{X_{ita}^r Y_{it}}{p_{it}^r} \middle |\ W\rb],\ \ho{E} \lb[\frac {X_{i't'a'}^r Y_{i't'}}{p_{i't'}^r} \middle|\ W\rb]\rb) \notag\\
&= \frac 1{p_{it}^r\cdot p_{i't'}^r}{\rm Cov}\lb(\ho{E} \lb[X_{ita}^r Y_{it}\middle|\ W_{{\cal N}^{t-\uptau^\star}(it)}\rb] + \eps,\ \ho{E}\lb[X_{i't'a'}^r Y_{i't'}\middle|\ W_{{\cal N}^{t'-\uptau^\star}(i't')}\rb] + \eps'\rb),
\end{align}
where  
\[\eps= \ho{E}[X_{ita}^r Y_{it}\mid W] -\ho{E} [X_{ita}^r Y_{it}\mid  W_{{\cal N}^{t-\uptau^\star}(it)}]\] and 
\[\eps' = \ho{E}[X_{i't'a'}^r Y_{i't'}\mid W] - \ho{E} [X_{i't'a'}^r Y_{i't'}\mid W_{{\cal N}^{t'-\uptau^\star}(i't')}].\]
To bound \cref{eqn:052425f}, we need two observations. 
First, {\new by the definition of $\uptau^\star$ and  our assumption that $t'-r \ge \uptau^\star$, we have 
\[W_{{\cal N}^{t-\uptau^\star}(it)} \ind W_{{\cal N}^{t'-\uptau^\star}(i't')}\] 
and consequently,
\begin{align}\label{eqn:052325b}
\ho{E} \lb[X_{ita}^r Y_{it}\mid W_{{\cal N}^{t-\uptau^\star}(it)}\rb] \ind \ho{E} \lb[X_{i't'a'}^r Y_{i't'}\mid W_{{\cal N}^{t'-\uptau^\star}(i't')}\rb].
\end{align} 
Furthermore, by \cref{prop:TVdecay}, we have
\begin{align}\label{eqn:052325a}
\|\eps\|_\infty\le e^{-|t-\uptau^\star| / t_{\rm mix}} 
\quad {\rm and}\quad 
\|\eps'\|_\infty \le e^{-|t'-\uptau^\star| / t_{\rm mix}}. 
\end{align} }
Thus, combining the Cauchy-Schwarz \ineq\ with \cref{eqn:052325a,eqn:052325b},  we obtain
\[\eqref{eqn:052425f}\le \frac {0 +  2e^{-|t'-\uptau^\star| / t_{\rm mix}} + \lb(e^{-|t'-\uptau^\star| / t_{\rm mix}}\rb)^2}{p_{it}^r\cdot p_{i't'}^r}
\le \frac {3e^{-|t'-\uptau^\star| / t_{\rm mix}}.}{p_{it}^r\cdot p_{i't'}^r}.\]
The desired \ineq\ follows by applying the above to all four combinations of $a,a'$.\qed

\end{document}